\documentclass[runningheads]{llncs}

\usepackage[T1]{fontenc}
\def\doi#1{\href{https://doi.org/\detokenize{#1}}{\url{https://doi.org/\detokenize{#1}}}}
\usepackage{graphicx}
\usepackage{listings}
\usepackage{hyperref}       % hyperlinks

\usepackage{amsmath,amssymb}%,amsthm}
\usepackage{algorithm, algorithmic}
\usepackage{pifont}
\usepackage{cases}
\usepackage{subcaption,graphicx}
\usepackage{svg}
\usepackage{stackengine}    % circled symbols
\newtheorem{assumption}[theorem]{Assumption}

\newcommand{\Max}{\max\limits}

\newcommand{\one}{\mathbf{1}}

\renewcommand{\le}{\leq}

\renewcommand{\hat}{\widehat}

\DeclareMathOperator{\spn}{span}
\DeclareMathOperator{\kernel}{Ker}

\DeclareMathOperator{\proj}{Proj}

\newcommand{\R}{\mathbb{R}}

\newcommand{\mA}{{\bf A}}

\newcommand{\mC}{{\bf C}}

\newcommand{\mW}{{\bf W}}

\newcommand{\sP}{{\mathcal{P}}}

\newcommand{\sX}{{\mathcal{X}}}
\newcommand{\sY}{{\mathcal{Y}}}

\newcommand{\bb}{{\bf b}}

\newcommand{\bd}{{\bf d}}

\newcommand{\bt}{{\bf t}}

\newcommand{\bx}{{\bf x}}
\newcommand{\by}{{\bf y}}
\newcommand{\bz}{{\bf z}}

\newcommand{\ds}{\displaystyle}
\newcommand{\norm}[1]{\left\| #1 \right\|}

\newcommand{\braces}[1]{\left\{ #1 \right\}}

\usepackage[colorinlistoftodos,bordercolor=orange,backgroundcolor=orange!20,linecolor=orange,textsize=scriptsize]{todonotes}

\newcommand{\ar}[1]{\textcolor{blue}{#1}}

\newcommand{\diag }[0]{\operatorname{diag }}

\usepackage{mathtools}
\newcommand{\defeq}{\vcentcolon=}

\newcommand{\tx}{{\tilde x}}
\newcommand{\ty}{{\tilde y}}
\newcommand{\tz}{{\tilde z}}
\newcommand{\tb}{{\tilde b}}
\newcommand{\td}{{\tilde d}}
\newcommand{\tm}{{\tilde m}}
\newcommand{\tn}{{\tilde n}}
\renewcommand{\th}{{\tilde h}}

\newcommand{\tbx}{{\tilde \bx}}
\newcommand{\tby}{{\tilde \by}}
\newcommand{\tbz}{{\tilde \bz}}

\newcommand{\ubx}{{\overline \bx}}
\newcommand{\uby}{{\overline \by}}

\newcommand{\tA}{{\tilde A}}
\newcommand{\tC}{{\tilde C}}

\newcommand{\tmW}{{\tilde \mW}}
\newcommand{\tsX}{{\tilde \sX}}
\newcommand{\tsY}{{\tilde \sY}}

\newcommand{\usX}{{\overline \sX}}
\newcommand{\usY}{{\overline \sY}}

\newcommand{\ihalf}{{i + \frac{1}{2}}}

\graphicspath{ {pics} }

\def\pd#1{{\color{black}#1}} %Pavel's comments
\def\ar#1{{\color{black}#1}} % Alexander R.'s comments
\def\dy#1{{\color{black}#1}} % Demyan's edits 

\begin{document}
\title{Decentralized convex optimization under affine
% separable
constraints for power systems control\thanks{The work of D. Yarmoshik was supported by the program “Leading Scientific Schools” (grant no. NSh-775.2022.1.1). The work of A. Rogozin and A. Gasnikov was supported by Russian Science Foundation (project No. 21-71- 30005).}}
\titlerunning{Decentralized convex optimization under affine constraints}
\author{Demyan Yarmoshik\inst{1}\orcidID{0000-0003-1912-1040} \and
Alexander Rogozin\inst{1}\orcidID{1111-2222-3333-4444}  \and
Oleg. O. Khamisov \inst{2}\orcidID{0000-0001-9015-9442}\and
Pavel Dvurechensky\inst{3}\orcidID{0000-0003-1201-2343} \and
Alexander Gasnikov\inst{1,4,5}\orcidID{0000-0002-7386-039X}}
\authorrunning{D. Yarmoshik et al.}
\institute{Moscow institute of physics and technology, Dolgoprudny, Institutskii ave., 9, Russia 
\email{\{yarmoshik.dv,aleksandr.rogozin,gasnikov.av\}@phystech.edu} \and Skolkovo Institute of Science and Technology, Moscow, Russia \email{Oleg.Khamisov@skolkovotech.ru}\and
Weierstrass Institute for Applied Analysis and Stochastics, Berlin, Germany \email{pavel.dvurechensky@wias-berlin.de}\and
Institute for Information Transmission Problems of the Russian Academy of Sciences (Kharkevich Institute), Russia \and 
Caucasus Mathematical Center, Adyghe State University, Maikop, Russia 
}

\maketitle              % typeset the header of the contribution
\begin{abstract}

Modern power systems are now in continuous process of massive changes. Increased penetration of distributed generation, usage of energy storage and controllable demand require introduction of a new control paradigm that does not rely on massive information exchange required by centralized approaches. Distributed algorithms can rely only on limited information from neighbours to obtain an optimal solution for various optimization problems, such as optimal power flow, unit commitment etc.

%One of the challenges in control of modern power systems is the task of solving an optimization problem over distributed network without sharing local information.

As a generalization of these problems we consider the problem of decentralized minimization of the smooth and convex partially separable function 
$f = \sum_{k=1}^l f^k(x^k,\tx)$
under the coupled 
$\sum_{k=1}^l (A^k x^k - b^k) \leq 0$
and the shared
$\tilde{A} \tilde{x}  - \tilde{b} \leq 0$
affine constraints, 
% $\sum_{k=1}^l (C^k x^k - d^k) \leq 0$,
where the information about $A^k$ and $b^k$ is only available for the $k$-th node of the computational network.
% This is a general formulation that covers aforementioned practical optimization problems.
%This problem covers multiple applications in energy systems control, e.g. direct current optimal power flow problem (DC-OPF) and schedule correction problem (??).

One way to handle the coupled constraints in a distributed manner is to rewrite them in a distributed-friendly form using the Laplace matrix of the communication graph and auxiliary variables (Khamisov, CDC, 2017).
% as described in \cite{Ol2,Ol4}.
Instead of using this method we reformulate the constrained optimization problem as a saddle point problem (SPP) and utilize the consensus constraint technique to make it distributed-friendly.
Then we provide a complexity analysis for state-of-the-art SPP solving algorithms applied to this SPP.
% We also show that 
% the method 
% of Khamisov 
% from \cite{Ol2,Ol4} 
% is equivalent to the standard consensus constraint technique from the SPP perspective. 

% We extend this approach to a broader class of problems which includes problems with shared variables (what makes the object function and the constraints inseparable on them).
% We applied 
% We apply saddle point approach for the convex case and show that the reformulation from [] is equivalent to the standard consensus constraint from the saddle point problem perspective.
% In the case of strongly convex object function and low-dimensional argument we utilize optimal dual methods.
% As an alternative for high-dimension strongly convex problems without the inequality-type constraints we apply an optimal algorithm from \cite{salim2021optimal}.
% All the proposed methods are supplied with the complexity estimations.

% The abstract should briefly summarize the contents of the paper in

% \keywords{First keyword  \and Second keyword \and Another keyword.}
\end{abstract}
\section{Introduction}

%Motivation
Optimal operation of power systems relies heavily on the ability of system operator to solve efficiently a number of optimization problems such as optimal power flow, unit commitment, as well as a number of online problems such as frequency and voltage control. Traditionally such problems were solved by System Operators in a centralized way. However, recent developments in implementation of distributed energy sources, storage systems and possibility of demand response can be effectively controlled by distributed algorithms. Such approach has a number of potential benefits, namely reduction of necessary communications between agents, increased robustness with respect to malfunction of any agent and possibility to increase cybersecurity and privacy of each agent. 

%Literature review
The detailed surveys on the application of distributed algorithms in power systems is given in \cite{molzahn2017survey,Patari2021}. 
\pd{These applications often lead to the necessity of solving an optimization problem, which can be formulated as distributed optimization problem with coupled constraints.}
Distributed approaches for optimization problems with coupled constraints can be separated into two main groups: (i) primal, dual or primal-dual consensus algorithms \cite{dc2,dc3,pdc1,pdc2,necoara2014distributed,NECOARA2015209,NECOARA2011756}; (ii) ADMM-based algorithms \cite{erseghe2014distributed,rostampour2019distributed,admm1,admm2}.

%Novelty

%Algorithm description
In this paper we propose a novel optimization approach for convex optimization problems with coupled linear equality and inequality constraints. Here introduction of specially placed Laplace matrices is used to model communications between neighboring agents in a computational network described as a connected graph. In the core of our approach lies: 1) the reduction of the decentralized optimization problem with constraints to decentralized saddle point problem; 2) applying decentralized Mirror Prox algorithm from \cite{rogozin2021decentralized} to solve the obtained saddle point problem. We obtain the same rate of convergence $\sim 1/N$ ($N$ -- number of communication steps / oracle calls) as the best known competitors, like ADMM \cite{lan2020first}. \pd{The main benefit of our approach is that the local optimization problem at each node is much simpler than in the ADMM-based approaches since we use only gradient oracle instead of complicated proximal mapping which may require a matrix inversion.
Compared to the dual algorithms of \cite{necoara2014distributed,NECOARA2015209,NECOARA2011756}, we consider a more general setting in which the objective may be non-separable and there are local linear constraints at each node of the computational network.
}
%But we use much cheaper gradient oracle, rather than proximal oracle in ADMM. 

% (??? need to describe solution algorithm ???)

\section{Problem Statement}
Let us consider the following optimization problem:
\begin{subequations}\label{eq:cp:p}
\begin{equation}
    %\min_{x\in\mathbb{R}^n}\frac{1}{2}x^\top Ax+b^\top x,\;A\succ0,
    \min_{x\in\mathbb{R}^n}f(\bx),
\end{equation}
\begin{equation}\label{eq:cp:p:ec}
    A'\bx - b'=0,\;A'\in\mathbb{R}^{m\times n},\;b'\in\mathbb{R}^m,
\end{equation}
\begin{equation}\label{eq:cp:p:ic}
    C'\bx - d'\le0,\;C'\in\mathbb{R}^{h\times n},\;d'\in\mathbb{R}^h,
\end{equation}
\end{subequations}
where $f:\mathbb{R}^n\rightarrow\mathbb{R}$ is a differentiable strictly convex function.
It is assumed that constraints \eqref{eq:cp:p:ec} and \eqref{eq:cp:p:ic} are consistent and there exists a unique solution $x^*$. Thus, Karush--Kuhn--Tucker (KKT) conditions are necessary and sufficient optimality conditions.

Let us now consider the case, when problem \eqref{eq:cp:p} must be solved by a multi-agent network with $l$ agents connected by a graph defined by a Laplacian matrix $W$. 
For this case, we assume, that  each agent seeks to find its own subvector $x^k\in\mathbb{R}^{n_k}$, $k\in\{1,\dots,l\}$ ($\sum_{k=1}^ln_k=n$) and the shared vector
$\tx \in \R^{\tn}$. We denote vector of private variables by $\bx = (x^{1\top},\dots,x^{l\top})^\top$.
Additionally, function $f$ is partially separable:
\begin{equation*}
    f(\bx, \tx) = \sum_{k=1}^lf^k(x^k, \tx)
\end{equation*}
and each $f^k$ is known only to agent $k$. Each agent has partial information $A^k\in\mathbb{R}^{m\times n_k}$, $b^k\in\mathbb{R}^m$, $C^k\in\mathbb{R}^{h\times n_k}$ and $d^k\in\mathbb{R}^h$ about constraints' parts corresponding only to variables $x^k$:
$A \defeq [A^1,\dots,A^l, ]$, $b \defeq \sum_{k=1}^lb^k$, $C \defeq [C^1,\dots,C^l]$ and  $n \defeq \sum_{k=1}^ln_k$. 
Additionally we assume that there are shared constraints with matrices 
$\tA \in \R^{\tm \times \tn}$,
$\tC \in \R^{\th \times \tn}$
and vectors 
$\tb \in \R^{\tm}$,
$\td \in \R^{\th} $
which are known to all agents. 

As a result, each agent $k$ has only its own part of the objective function $f^k(x^k,\tilde{x})$ and parts of the coupled equality and inequality constraints respectively: $A^kx^k-b^k$ and $C^kx^k-d^k$.

Therefore, we have an optimization problem of the following form:

\begin{subequations}\label{eq:cpc:p}
\begin{align}
	\min_{x\in\R^{n+\tn}}~ &\sum_{k=1}^l f^k(x^k,\tx), \\
	\text{s.t. }&\sum_{k=1}^l (A^k x^k - b^k) = 0, \\
	&\sum_{k=1}^l (C^k x^k - d^k) \leq 0,\\
	&\tA\tx-\tb = 0,\label{eq:cpc:p:sheq}\\
	&\tC\tx-\td\le0\label{eq:cpc:p:shin}.
\end{align}
\end{subequations}
Here $\tx\in\mathbb{R}^{\tn}$ is a subvector of $x$ that contains global variables used by all agents. 

% In order to modify this problem for further distributed solution let us introduce approximate copies of variables $\tx$ for each agent: $\tx^k\in\mathbb{R}^{\tn},\;k\in\{1,\dots,l\}$.

\section{Mathematical setting}

\begin{assumption}\label{ass:fk}
For every $k=1,\ldots,l$
\begin{enumerate}
    \item $f^k(x^k, \tx)$ is differentiable.
    \item (Convexity) $ \forall x^k, x'^k \in \sX^k, \forall \tx, \tx' \in \tilde \sX $
    $$f^k(x'^k, \tx') \geq f^k(x^k, \tx) + \left\langle \nabla f^k(x^k, \tx), 
    \begin{pmatrix}
    x'^k - x^k\\ \tx' - \tx
    \end{pmatrix}%(x^k, \tx)^\top
    \right\rangle. $$
    \item (Lipschitz smoothness) 
    $$\left\|\nabla f^k(x'^k, \tx') - \nabla f^k(x^k, \tx)\right\| \leq L_k \left\|
    \begin{pmatrix}
     x'^k- x^k \\ \tx' - \tx
    \end{pmatrix}
    \right\|.$$
\end{enumerate}
\end{assumption}

\begin{assumption}\label{ass:blocks}
Variable $x$ is subject to block constraints: 
$x^k \in \prod_{i=1}^{n_k} [\xi^{k,i}, \eta^{k,i}] = \sX^k$, $\xi^{k,i}, \eta^{k,i} \in \R$ 
and 
$\tx \in \prod_{i=1}^{\tn} [\tilde \xi^{i}, \tilde \eta^{i}] = \tilde\sX$, $\tilde \xi^{i}, \tilde \eta^{i} \in \R$ .
\end{assumption}
This is a natural assumption since in a real-world system maximal and minimal values of every control and auxiliary variable are limited.
Let us also denote
\begin{itemize}
    \item  $\lambda_{max}(A), \lambda_{min}^+(A)$~--- the largest and the smallest positive eigenvalues of a matrix A.
    \item $\sigma_{max}(A) = \sqrt{\lambda_{max}(A^\top A)}$ and $\sigma^+_{min}(A) = \sqrt{\lambda^+_{min}(A^\top A)}$~--- the largest and the smallest positive  singular value of a matrix A.
    \item $\chi(A) = \frac{\sigma_{max}(A)}{\sigma^+_{min}(A)}$~--- condition number of a matrix $A$ on $\left(\kernel{A}\right)^\top$.
    \item $\proj_S(x)$~--- projection of $x$ onto a set $S$.
\end{itemize}

The key instrument in separating shared variables and coupled constraints is introducing the consensus constraint with the help of matrix $W$ defined as follows:
\begin{enumerate}
    \item ${ W}$ is symmetric positive semi-definite matrix.
	\item (Network compatibility) For all $i,j = 1,\dots,l$ the entry of $ W$: $[{ W}]_{ij} = 0$ if $i\ne j$ and there is no edge in the communication graph between nodes $i$ and $j$. This property allows to perform multiplications by $W$ in a distributed manner (only using information from neighbours in the communication graph).
	\item (Kernel property) For any $v = [v_1,\ldots,v_m]^\top\in\R^m$, ${ W} v = 0$ if and only if $v_1 = \ldots = v_m$, i.e.  $\kernel { W} = \spn\braces{\one}$. This property allows to rewrite pairwise equality constraint in a distributed way.
\end{enumerate}

An example of matrix satisfying this assumption is the graph Laplacian ${ W}{\in \mathbb{R}^{m\times m}}$:
{\begin{align*}
[W]_{ij} \triangleq \begin{cases}
-1,  & \text{if } (i,j) \in E,\\
\text{deg}(i), &\text{if } i= j, \\
0,  & \text{otherwise,}
\end{cases}
\end{align*}}
where $\text{deg}(i)$ is the degree of the node $i$, i.e., the number of neighbors of the node. 

Matrix $W$ can be used to rewrite pairwise equality of scalars.
To rewrite pairwise equality of vector variables with equal dimesion we will use the following extension of matrix $W$, called \textit{communication matrix}:
\begin{equation}\label{eq:boldsymW}
    \mW  = W \otimes I_d,
\end{equation}
where $\otimes$ denotes the Kronecker product and $d$ is the dimension of the vector variables.

\section{Distributed saddle point problem formulation}

\subsection{Saddle point problem and consensus constraints}
% \subsection{Centralized case}

% \cite{rogozin2021decentralized}
We reformulate problem \eqref{eq:cp:p} as saddle point problem:
\begin{multline}\label{eq:sf-sp}
    \min_{\bx, \tx} \max_{\substack{\lambda, \tilde \lambda\\ \mu, \tilde \mu \geq 0}} \sum_{k=1}^l 
    \left[ f^k(x^k, \tx) + 
    \lambda^\top (A^k x^k - b^k) + \mu^\top (C^k x^k - d^k) \right]  +\\
    \tilde \lambda^\top (\tA\tx -\tb) +
    \tilde \mu^\top (\tC \tx - \td) 
    .
\end{multline}

Let us unify the analysis of equality and inequality constraints by stacking Lagrange multipliers
$\lambda$ and $\mu$ in a single dual variable
\begin{equation*}
y = \begin{pmatrix} \lambda \\ \mu \end{pmatrix},~
y \in \sY = \R^m \times \R^h_+ .
\end{equation*}

And similarly we introduce the joined dual variable for the coupled constraints:

\begin{equation*}
\ty = \begin{pmatrix} \tilde \lambda \\ \tilde \mu \end{pmatrix},~
\ty \in \tsY = \R^{\tm} \times \R^{\th}_+ .
\end{equation*}

% \subsection{Distributed reformulation}
To solve this saddle point problem in a distributed manner we have to separate dual variables $y$
by making their copies at each node and introducing consensus constraint into the saddle point problem, as described in  \cite{rogozin2021decentralized}.
That brings us to the following formulation:

\begin{equation}%\label{eq:sf-sp:d}
    \min_{\bx, \tx, \bz} \max_{\by, \ty}  \sum_{k=1}^l \left[f^k(x^k, \tx) + 
    y^{k\top} 
    \begin{pmatrix} A^k x^k - b^k  \\C^k x^k - d^k \end{pmatrix}  
    \right] + 
    \bz^\top \mW \by +
    \ty^{\top} 
    \begin{pmatrix} \tA \tx - \tb  \\
    \tC \tx - \td \end{pmatrix}  
    % = \min_{x,z} \max_y~  \bz^\top \mW\by + \sum_{k=1}^l g^k(x^k, y^k).
\end{equation}

To separate the terms corresponding to the shared constraints \eqref{eq:cpc:p:sheq}, \eqref{eq:cpc:p:shin} we should go back to the optimization problem \eqref{eq:cpc:p} and do the same trick with them:
make a copy of $\tx$ at each node and introduce consensus constraint.
So we transform \eqref{eq:cpc:p:sheq}, \eqref{eq:cpc:p:shin} into equivalent system
\begin{subequations}\label{eq:cpc:p:W}
\begin{align}
	&\tA\tx^k +\tb = 0, \; k \in \{1, \ldots, l \} ,\label{eq:cpc:p:W:sheq}\\
	&\tC\tx^k +\td\le0, \; k \in \{1, \ldots, l \}\label{eq:cpc:p:W:shin},\\
	&\tmW \tbx = 0,
\end{align}
\end{subequations}
where $ \tbx = (\tx^{1\top},\dots,\tx^{l\top})^\top, \; \tmW = W \otimes I_\tn$.

Note, that each node can handle constraints \eqref{eq:cpc:p:W:sheq} and \eqref{eq:cpc:p:W:shin} independently, so we don't have to introduce additional consensus constraints over corresponding dual variables in the final saddle point problem:

\begin{multline}
\label{eq:sf-sp:d}
    \min_{\bx, \tbx, \bz} \max_{\by, \tby, \tbz}
    \sum_{k=1}^l \left[f^k(x^k, \tx^k) + 
    y^{k\top} 
    \begin{pmatrix} A^k x^k - b^k  \\C^k x^k - d^k \end{pmatrix}  +
    \ty^{k\top} 
    \begin{pmatrix} \tA \tx^k - \tb  \\
    \tC \tx^k - \td \end{pmatrix}  
    \right] + 
    \bz^\top \mW \by +
    \tbz^\top \tmW \tbx \\
    % = \min_{x,z} \max_y~  \bz^\top \mW\by 
    = \min_{\bx, \tbx, \bz} \max_{\by, \tby, \tbz}
    \sum_{k=1}^l g^k(x^k, \tx^k, y^k, \ty^k) +  \bz^\top \mW \by +
    \tbz^\top \tmW \tbx,
\end{multline}
%  x -> x
%  y -> r
%  z -> u
%  g -> f

where $\by = (y^{1\top}, \ldots, y^{l\top})^\top$, $\mW = W \otimes I_{m+h}$.

\dy{We will also use the following notation:
\begin{equation}
    G(\bx, \tbx, \by, \tby) = \sum_{k=1}^l g^k(x^k,\tx^k, y^k, \ty^k),
\end{equation}
and
\begin{equation}
    G_w(\bx, \tbx, \by, \tby) = G(\bx, \tbx, \by, \tby)  + 
    \bz^\top \mW \by +
    \tbz^\top \tmW \tbx.
\end{equation}
}

\subsection{Comparison with \cite{Ol2}, \cite{Ol4} }
In this subsection we show the equivalence of our approach and approach from \cite{Ol2}, \cite{Ol4} from the perspective of saddle point problems.
Since the shared variables $\tx$ are handled in the same way in both approaches (by introducing the constraint $\tmW \tbx = 0$ into the optimization problem), we consider the case without shared variables and only with equality-type constraints to simplify the derivations.

Let us introduce a set of new matrices and vectors:
\begin{equation}\label{eq:zero}
\begin{gathered}
    \mA = \diag(A^1,\dots,A^l), \bb = (b^{1\top},\dots,b^{l\top})^\top,\\
    % \mC = \diag(C^1,\dots,C^l), \bd = (d^{1\top},\dots,d^{l\top})^\top,\\
\end{gathered}
\end{equation}
\begin{equation}\label{eq:g}
\begin{gathered}
    W^{mk} = \diag(W_{k\bullet},\dots,W_{k\bullet})\in\mathbb{R}^{m\times ml}, \mW^m = \left[\begin{array}{c}
        W^{m1}\\
        \vdots\\
        W^{ml}
    \end{array}\right]\in\mathbb{R}^{ml\times ml},\\
    % W^{hk} = \diag(W_{k\bullet},\dots,W_{k\bullet})\in\mathbb{R}^{h\times hl}, \mW^h = \left[\begin{array}{c}
    %     W^{h1}\\
    %     \vdots\\
    %     W^{hl}
    % \end{array}\right]\in\mathbb{R}^{hl\times hl}.
\end{gathered}
\end{equation}

In \cite{Ol2}, \cite{Ol4} the following distributed-friendly reformulation of problem \eqref{eq:cpc:p} is proposed, and its equivalence to the original problem is shown:
%$W^m\in\mathbb{R}^{ml\times l}$ and $W^h\in\mathbb{R}^{hl\times l}$ which a constructed by repeating each row of $W$ $m$ and $h$ times respectively. 
% Then problem \eqref{eq:dp:p} can be represented in matrix form:
\begin{subequations}\label{eq:dp:mp}
\begin{equation}
    %\min_{x\in\mathbb{R}^n,y\in\mathbb{R}^{l\times m},z\in\mathbb{R}^{l\times h}}\frac{1}{2}x^\top Ax+b^\top x,
   \min_{\bx\in\mathbb{R}^n,\by\in\mathbb{R}^{ml}
    % ,z\in\mathbb{R}^{hl}
    } \left\{f(\bx) =  \sum_{k=1}^l f^k(x^k) \right\}
     ,
\end{equation}
\begin{equation}\label{eq:dp:mp:e}
    \mA \bx - \bb + \mW^m \by = 0.
\end{equation}
% \begin{equation}\label{eq:dp:mp:i}
%     \mC x + \bd + \mW^h z\le0.
% \end{equation}
\end{subequations}
Here a sort of consensus constraint is integrated directly into the minimization problem, which differs from our technique of adding consensus constraint into the corresponding saddle point problem.
Note also that $\mW^m$ and $\mW$ differ in their structure (the way of constructing matrix $W$ for using it with multi-dimensional variables).

% Let us write Lagrangian for distributed-friendly reformulation of problem \eqref{eq:cp:p} given by problem \eqref{eq:dp:p} (problem \eqref{eq:dp:mp})
% and then solve saddle point problem without explicitly introducing additional constraints of type $\mW^m\bz = 0$ 
% % but it follows from $\nabla_y L = 0$ 
% (for simplicity we skip inequality constraints):
The saddle point problem corresponding to the minimization problem \eqref{eq:dp:mp} is

\begin{equation}\label{eq:rf-lg-uf}
    \min_{\bx,\by} \max_{\bz} L(\bx, \by, \bz) =
    \min_{\bx,\by} \max_{\bz} f(\bx) + 
     \bz^\top \left(\mA\bx - \bb + \mW^m \by \right).
\end{equation}

Let us now compare this problem with the saddle point problem \eqref{eq:sf-sp:d}.
By rewriting sum in \eqref{eq:sf-sp:d} and using the symmetry of $\mW$ we have
\begin{align}\label{eq:sf-sp:d:m}
    \min_{\bx,\bz} &\max_\by  \sum_{k=1}^l \left[f^k(x^k) + y^{k\top} (A^k x^k - b^k)  \right] + \bz^\top \mW \by \notag \\
    &= \min_{\bx,\bz} \max_{\by} f(\bx) + \by^{\top} (\mA\bx - \bb) + \bz^\top \mW \by \notag \\
    &= \min_{\bx,\bz} \max_{\by} f(\bx) + \by^{\top} (\mA\bx - \bb + \mW\bz).
\end{align}
Since $\mW$ and $\mW^m$ differ only in the arrangement of columns,
problems \eqref{eq:rf-lg-uf} and \eqref{eq:sf-sp:d:m} differ only in the arrangement of components of maximized variables.
Therefore, both approaches leads to the same saddle point problem.

% Reformulated optimization problem for common variables case \eqref{e:dp_full} also leads to the same saddle point problem as straightforward approach \eqref{eq:sf_l_full}, \eqref{eq:sf_sp_full}.
% It is so because constraints on common variables $\tilde \bx$ introduced in the same way in both approaches, and terms for local variables constraints differ only in arrangement of components of maximized variables, as shown in \eqref{eq:sf-sp:d:m}. 

\section{Algorithm}
We use classical Extragradient algorithm from \cite{korpelevich1976extragradient}.
Being applied to the problem \eqref{eq:sf-sp:d} it converges to the solutions of the primal and the dual problems as will be shown in the next sections.
Here we describe it in an explicit form, so it is ready to be applied to the problem \eqref{eq:cpc:p}, see Algorithm \ref{alg:extragrad}.

Note, that the projection in our case is a simple clipping and can be performed independently for each component of the variable.

\begin{algorithm}[ht!]
\caption{Decentralized Extragradient for problem \eqref{eq:cpc:p}}
{
\begin{algorithmic}[1]\label{alg:extragrad}
    \STATE {Initialize $\bx_0 \in \sX, \by_0 = \bz_0 = \Vec{0}_{l(m+h)}, \tbx_0 \in \tsX^l, \tby = \Vec{0}_{l(\tm+\th)}, \tbz_0 = \Vec{0}_{l \tn}$}
    \FOR{$i = 0, \ldots, N - 1$}
        \STATE{Compute $\bz'_i = \mW \bz_i$, $\by'_i = \mW \by_i$, $\tbz'_i = \tmW \tbz_i$, $\tbx'_i = \tmW \tbx_i$}.
        \STATE{Make intermediate gradient step
        \begin{align*}
            &x_{\ihalf}^k = \proj_\sX \left( x^k_i - h \nabla_{x^k} f^k(x_i^k, \tx_i^k) - h (A^{k\top}, C^{k\top})y^k_{i} \right)\\
            &\tx_\ihalf^k = \proj_\tsX \left(\tx^k_i - h \nabla_{\tx^k} f^k(x_i^k, \tx_i^k) - h \tz'^k_i \right)\\
            &y_\ihalf^k = \proj_\sY \left(
            y^k_i + h \begin{pmatrix} A^k x_i^k - b^k  \\C^k x_i^k - d^k \end{pmatrix} +h z'^k_i \right)
            \\
            &\ty_\ihalf^k =\proj_\tsY \left( 
            \ty^k_i + h \begin{pmatrix} \tA^k \tx_i^k - \tb^k  \\\tC^k \tx_i^k - \td^k \end{pmatrix} \right)
            \\
            &z_\ihalf^k = z^k_i - h y'^k_{i}\\
            &\tz_\ihalf^k = \tz^k_i + h \tx'^k_i \\
        \end{align*}
        }
    
    \STATE{Compute $\bz'_\ihalf = \mW \bz_\ihalf$, $\by'_\ihalf = \mW \by_\ihalf$, $\tbz'_\ihalf = \tmW \tbz_\ihalf$, $\tbx'_\ihalf = \tmW \tbx_\ihalf$}.
    \STATE{Make gradient step
    \begin{align*}
        &x_{i + 1}^k = \proj_\sX \left( x^k_i - h \nabla_{x^k} f^k(x_\ihalf^k, \tx_\ihalf^k) - h (A^{k\top}, C^{k\top})y^k_\ihalf \right)\\
        &\tx_{i + 1}^k = \proj_\tsX \left( \tx^k_i - h \nabla_{\tx^k} f^k(x_\ihalf^k, \tx_\ihalf^k) - h \tz'^k_\ihalf \right)\\
        &y_{i + 1}^k = \proj_\sY \left(
        y^k_i + h \begin{pmatrix} A^k x_\ihalf^k - b^k  \\C^k x_\ihalf^k - d^k \end{pmatrix} +h z'^k_\ihalf \right)
        \\
        &\ty_{i + 1}^k =\proj_\tsY \left( 
        \ty^k_i + h \begin{pmatrix} \tA^k \tx_\ihalf^k - \tb^k  \\\tC^k \tx_\ihalf^k - \td^k \end{pmatrix} \right)
        \\
        &z_{i + 1}^k = z^k_i - h y'^k_\ihalf\\
        &\tz_{i + 1}^k = \tz^k_i + h \tx'^k_\ihalf \\
    \end{align*}
    }
        
    %     \STATE{Each node computes
    %         \begin{align*}
    %             &x_i^{k+\frac{1}{2}} = {\rm Mirr}\left(
    %             \nabla_{x} f_i^k
    %             +{\gamma_x}\tilde z_i^k; x_i^k; \bar\sX\right),~~~~
    %             s_i^{k+\frac{1}{2}} = {\rm Mirr} \cbraces{{\gamma_y}\tilde y_i^k; s_i^k; \R^{d_y}}, \\
    %             &y_i^{k+\frac{1}{2}} = {\rm Mirr} \cbraces{-\nabla_r f_i^k -  {\gamma_y} \tilde s_i^k; y_i^k; \bar\sY},~
    %             z_i^{k+\frac{1}{2}} = {\rm Mirr} \cbraces{-{\gamma_x} \tilde x_i^k; z_i^k; \R^{d_x}}.
    %         \end{align*}
    %         }
    %     \STATE{Compute $\tilde\bx^{k+\frac{1}{2}} = \mW_x\bx^{k+\frac{1}{2}}$;~ 
    %     % \tilde\bs^{k+\frac{1}{2}} = \mWy\bs^{k+\frac{1}{2}};~ 
    %     % \tilde\by^{k+\frac{1}{2}} = \mWy\by^{k+\frac{1}{2}};~ 
    %     % \tilde\bz^{k+\frac{1}{2}} = \mWx\bz^{k+\frac{1}{2}}$.
    %     }
    %     \STATE{Each node computes
    %         \begin{align*}
    %             &x_i^{k+1} = {\rm Mirr}\left(
    %             \nabla_{p} f_i^{k+\frac{1}{2}}
    %             +{\gamma_x}\tilde z_i^k; x_i^k; \bar\sX\right),~~~~
    %             s_i^{k+1} = {\rm Mirr} \cbraces{{\gamma_y}\tilde y_i^{k+\frac{1}{2}}; s_i^k; \R^{d_y}}, \\
    %             &y_i^{k+1} = {\rm Mirr} \cbraces{-\nabla_r f_i^{k+\frac{1}{2}} -  {\gamma_y} \tilde s_i^k; y_i^k; \bar\sY},~
    %             z_i^{k+1} = {\rm Mirr} \cbraces{-{\gamma_x} \tilde x_i^{k+\frac{1}{2}}; z_i^k; \R^{d_x}}.
    %         \end{align*}
    %         }
    \ENDFOR
    \ENSURE{For $\bt\in\braces{\bx, \by, \bz, \tbx, \tby, \tbz}$ compute $\ds\hat \bt^{N} = \frac{1}{N}\sum_{k=0}^{N-1} \bt^{k+\frac{1}{2}}$.}
\end{algorithmic}
}
\end{algorithm}

\section{Smoothness and domain size analysis}
In this section we will perform some technical analysis to obtain the relations between parameters of the input data to the problem (object functions and constrains) and parameters of Extragradient's convergence rate.

\subsection{Bounds on $\|y^*\|$, $\| \ty^*\|$}
To calculate Lipschitz smoothness constants of the problem we have to localize $y^*$ (dual part of solution of the initial saddle problem \eqref{eq:sf-sp}, which is also a solution to the dual problem under our assumptions),
i.e. find $R_y$ such that $\sY$ lies in a ball in $\R^{m}$ with center in $0$ and radius $R_y$, and $y^* \in \sY$. 
% $\sup_{y^0, y^{0'} \in \sY^0} \|y - y'\| \leq R_y$.
From optimality conditions for dual problem of  \eqref{eq:cpc:p} 
\begin{equation}\label{eq:zerograd_1}
\nabla_\bx L = \nabla_\bx f + (A^\top, C^\top) y^* = 0 
\end{equation}
\begin{equation}\label{eq:zerograd_2}
\nabla_{\tx} L = \nabla_{\tx} f + (\tA^\top, \tC^\top) \ty^*  = 0.
\end{equation}
Since for any 
$y \in \ker{A^T}$ vector 
$y^* + y$ is also a solution, we consider only solution with the smallest norm (it's enough for saddle point problem solution's quality criteria and convergence analysis), i.~e. 
$y^* \in \left(\ker{A^T}\right)^\bot$. 

Therefore

% We obtain $R_y$ by considering following inequality \cite{rogozin2021decentralized}:
\begin{equation*}
   \|y^*\|^2 \leq \frac{\|\nabla_\bx f(\bx^*, \tx^*)\|^2}
   {(\sigma_{min}^+\left((A^\top, C^\top)\right))^2},
\end{equation*}
\begin{equation*}
   \|\ty^*\|^2 \leq \frac{\|\nabla_{\tx} f(\bx^*, \tx^*)\|^2}
   {(\sigma_{min}^+((\tA^\top, \tC^\top)))^2},
\end{equation*}
where $\sigma_{min}^+(A) = 
\sqrt{ \min\{ \lambda > 0 : \exists x \neq 0 : AA^T x = \lambda x \} }$.
Hence we get  

\begin{lemma}
Saddle point problem \eqref{eq:sf-sp:d}, which is unconstrained on variables $\by, \tby$, is equivalent to the same problem with constraints $\|\by\| \leq R_\by$ and $\|\tby \| \leq R_{\tby}$, where

\begin{equation}\label{eq:ry-sf}
	R_\by = \sqrt{l}~\frac{\Max_{\bx \in \sX, \tx \in \tsX } \|\nabla_\bx f(\bx, \tx)\|}
   {\sigma_{min}^+\left((A^\top, C^\top)\right)},
	R_{\tby} = \sqrt{l}~\frac{\Max_{\bx \in \sX, \tx \in \tsX } \|\nabla_{\tx} f(\bx, \tx)\|}
	{\sigma_{min}^+\left((\tA^\top, \tC^\top)\right)}.
\end{equation} 
\end{lemma}

\subsection{Bounds on $\|z^*\|$, $\| \tz^*\|$}
Next we want to find constants for Euclidean-case bounds for Theorem 3.5 \cite{rogozin2021decentralized}.
To specify, how the convergence rate depends on problem's parameters, we need to find scalars $M_y, M_{\tx}, L_{xx}, L_{yx}, L_{xy} L_{yy}$, determined by inequalities

\begin{subequations}
\begin{align}
	&\| \nabla_y g^k(x^k, \tx^k, y^k, \ty^k) \| \leq M_y ~\forall k, x_k \in \sX_k, y_k \in \sY \\
	&\| \nabla_{\tx} g_k(x^k, \tx^k, y^k, \ty^k) \| \leq M_{\tx} ~\forall k, x_k \in \sX_k, y_k \in \sY \\
	&\|\nabla_x G(\ubx, \uby) - \nabla_x G(\ubx', \uby)\| \leq L_{xx}\|\ubx - \ubx'\| ~\forall \ubx, \ubx' \in \usX, \uby \in \usY \\
	&\|\nabla_x G(\ubx, \uby) - \nabla_x G(\ubx, \uby')\| \leq L_{xy}\|\uby - \uby'\| ~\forall \ubx \in \usX, \uby, \uby' \in \usY \\
	&\|\nabla_y G(\ubx, \uby) - \nabla_y G(\ubx', \uby)\| \leq L_{yx}\|\ubx - \ubx'\| ~\forall \ubx, \ubx' \in \usX,~ \forall  \uby \in \usY \\
	&\|\nabla_y G(\ubx, \uby) - \nabla_y G(\ubx, \uby')\| \leq L_{yy}\|\uby - \uby'\| ~\forall \ubx \in \usX,~ \forall \uby, \uby' \in \usY, 
\end{align}
\end{subequations}
where $\ubx = (\bx^\top, \tbx^\top)^\top$ and $\uby = (\by^\top, \tby^\top)^\top$

% Let us denote 
% $\sigma_{max}(A) = \sqrt{\lambda_{max}(A^\top A)}$~--- the largest singular value of A.
% Then
\newcommand{\kol}{{k\in \{1, \ldots, l \}}}
\newcommand{\iom}{{i\in \{1, \ldots, m \}}}
% \begin{equation*}

\dy{By using the triangle inequality}
%M_y
\begin{multline*}
	\| \nabla_{y^k} g^k(x^k, \tx^k, y^k, \ty^k) \| = \|\nabla_{y^k} y^\top (A^k x^k - b^k)\| =
% 	\|A^k x^k + b^k\|
    \left\|\begin{pmatrix} A^k x^k - b^k  \\C^k x^k - d^k \end{pmatrix} \right\|
	\leq  \\
	\max_{\kol} \left\{
% 	R^k_x = \max_{x^k \in \sX^k}\|x^k\|
    \sigma_{max}\left(({A^k}^\top, {C^k}^\top)\right)R_{x^k} + \left\|({b^k}^\top, {d^k}^\top)\right\| 
	\right\}
	%\leq \max_{k\in \{1, \ldots, l \}, x^k \in \sX^i} \|A^k x^k + b^k\| =
	= M_y,
\end{multline*}
% \end{equation*}
\dy{and}
% M_\tx
\begin{multline*}
	\| \nabla_{\tx^k} g^k(x^k, \tx^k, y^k, \ty^k) \| =
    \| \begin{pmatrix}{\tA}^\top, {\tC}^\top \end{pmatrix}
    \ty^k \| 
	\leq   \sigma_{max}\left(\begin{pmatrix}{\tA}^\top, {\tC}^\top\; \end{pmatrix}\right)R_{\ty}   \\
% 	{(\sigma_{min}^+(A))^2}
	= \chi\left(({\tA}^\top, {\tC}^\top)\right)  \Max_{\bx \in \sX, \tx \in \tsX} \|\nabla_{\tx} f(\bx, \tx)\|
% 	\right\}
	= M_{\tx}.
\end{multline*}

Then by directly applying Lemma 4.2 in \cite{rogozin2021decentralized} we have

\begin{lemma}
Saddle point problem \eqref{eq:sf-sp:d}, which is unconstrained on variables $\bz, \tbz$, is equivalent to the same problem with constraints $\|\bz\| \leq R_\bz$ and $\|\tbz \| \leq R_{\tbz}$, where
\begin{equation}
    R_\bz = \frac{\sqrt{2l}M_y}{\lambda_{min}^+(\mW)},
    R_{\tbz} = \frac{\sqrt{2l}M_{\tbx}}{\lambda_{min}^+(\tmW)}.
\end{equation}
\end{lemma}

\subsection{Smoothness constants}
Let us find smoothness constants of function  $G$. From \eqref{eq:sf-sp:d} we have
\begin{equation*}
\nabla_x G(\ubx, \uby) - \nabla_x G(\ubx', \uby) = 
\begin{pmatrix}
\nabla f^1(x^1, \tx^1) - \nabla f^1(x^{1'}, \tx^{1'})\\
\vdots \\
\nabla f^l(x^l, \tx^l) - \nabla f^l(x^{l'}, \tx^{l'})\\
\end{pmatrix}.
\end{equation*}
By Assumption \ref{ass:fk}
\begin{equation*}
\begin{aligned}
&\|\nabla_x G(\ubx, \uby) - \nabla_x G(\ubx', \uby)\|^2 =
\sum_{k=1}^l \| \nabla f^k(x^k, \tx^k) - \nabla f^k(x^{k'}, \tx^{k'}) \|^2 \\
&\leq \sum_{k=1}^l L_k^2\left\|
    \begin{pmatrix}
     x'^k- x^k \\ \tx' - \tx
    \end{pmatrix}
    \right\|^2% \|x^k - x^{k'} \|^2
\leq \max_k L_k^2 \|\ubx - \ubx'\|^2.
\end{aligned}
\end{equation*}

Taking square root from both parts of the inequality we get 
\begin{equation*}
    L_{xx} = \max_\kol  L_k .
\end{equation*}

Similarly, for other variables
\begin{equation*}
\nabla_x G(\ubx, \uby) - \nabla_x G(\ubx, \uby')= 
\begin{pmatrix}
    ( {A}^\top, {C}^\top)\;
    (y - y') \\
    \begin{pmatrix}{\tA}^\top, {\tC}^\top\; \end{pmatrix}
    (\ty^{l} - \ty^{l'})  \\
    % % \begin{pmatrix}
    % ( {A^1}^\top, {C^1}^\top)\;
    % % \end{pmatrix}
    % (y^{1} - y^{1'}) \\ \vdots \\
    % % \begin{pmatrix}
    % ({A^l}^\top, {C^l}^\top) \;
    % % \end{pmatrix}
    % (y^{l} - y^{l'})  \\
    
    % \begin{pmatrix}{\tA}^\top, {\tC}^\top\; \end{pmatrix}
    % (\ty^{1} - \ty^{1'}) \\ \vdots \\
    % \begin{pmatrix}{\tA}^\top, {\tC}^\top\; \end{pmatrix}
    % (\ty^{l} - \ty^{l'})  \\
\end{pmatrix},
\end{equation*}
\begin{equation*}
\nabla_y G(\ubx, \uby) - \nabla_y G(\ubx', \uby)= 
\begin{pmatrix} 
    ( {A^1}^\top, {C^1}^\top)^\top (x^{1} - x^{1'}) \\
    % C^{1}(x^{1} - x^{1'}) \\
    \vdots \\
    ( {A^l}^\top, {C^l}^\top)^\top (x^{l} - x^{l'}) \\
    % C^{l}(x^{l} - x^{l'}) \\
   % 
    \begin{pmatrix}{\tA}^\top, {\tC}^\top\; \end{pmatrix}^\top
    (\tx^{1} - \tx^{1'}) \\
    \vdots \\
    \begin{pmatrix}{\tA}^\top, {\tC}^\top\; \end{pmatrix}^\top
    (\tx^{l} - \tx^{l'})
\end{pmatrix}.
\end{equation*}
and

\begin{equation*}
L_{xy} = 
\max \left\{\max_\kol \sigma_{max}\left( ({A^k}^\top, {C^k}^\top)\right), 
\; \sigma_{max} \left((\tA^\top, \tC^\top )\right) \right\} = L_{yx} ,
% \max_\kol \sigma_{max} (A^k)
\end{equation*}

% \begin{equation*}
% \begin{aligned}
% &\|\nabla_x G(x^0, y^0) - \nabla_x G(x^{0}, y^{0'})\| =
% \|\sum_{k=1}^l A^{kT}(y^{k} - y^{k'}) \| \leq \max_k \sigma_{max} (A^k) \|y^0 - y^{0'}\| = L_{xy} \|y^0 - y^{0'}\|
% \end{aligned}
% \end{equation*}

% \begin{equation*}
% \begin{aligned}
% &\|\nabla_y G(x^{0}, y^0) - \nabla_y G(x^{0'}, y^{0})\| =
% \|\sum_{k=1}^l A^{k}(x^{k} - x^{k'}) \| \leq \max_k \sigma_{max} (A^k) \|x^0 - x^{0'}\| = L_{yx} \|x^0 - x^{0'}\|
% \end{aligned}
% \end{equation*}

\begin{equation*}
\begin{aligned}
L_{yy} = 0.
\end{aligned}
\end{equation*}

\section{Main result}
Let us denote 
\ar{
\begin{align*}
L_\zeta =  2\cdot \text{max}\{
&R^2_{\bx\tbx} L_{x \tx, x \tx}, \;
R^2_{\by\tby} L_{y \ty, y \ty}, \; \\
&\sqrt{2} R_{\bx\tbx}R_{\by\tby}L_{x \tx, y \ty} + 2M_{\bx\tbx}R_{\bx\tbx} \frac{\lambda_{\max}(\tilde\mW)}{\lambda_{\min}^+(\tilde\mW)} + 2M_{\by\tby}R_{\by\tby} \frac{\lambda_{\max}(\mW)}{\lambda_{\min}^+(\mW)}\}.
\end{align*}
}
% \begin{equation}
% L_\zeta =  2 \max\{
% % (R^2_\sX +  R^2_{\tsX})L_{x \tx, x \tx},
% % (R^2_\sY +  R^2_{\tsY})L_{y \ty, y \ty}, 
% R^2_{\bx\tbx}L_{x \tx, x \tx}, \;
% R^2_{\by\tby}L_{y \ty, y \ty}, \;
% R_{\bx\tbx}R_{\by\tby}L_{x \tx, y \ty}, \;
% R_\by R_\bz \|\mW \|, \;
% R_{\tbx} R_{\tbz} \|\tmW \|
% \}.
% \end{equation}

Then, \ar{following the arguments presented in} Theorem 3.5 from \cite{rogozin2021decentralized}, \ar{we introduce $\zeta = (\bx^\top, \tbx^\top, \by^\top, \tby^\top, \bz^\top, \bz^\top)^\top$}. \ar{We also define a norm for $\zeta$ as follows:
\begin{align*}
    \norm{\zeta}^2 = \frac{\norm{\bx}^2}{R_\bx^2} + \frac{\norm{\tbx}^2}{R_\tbx^2} + \frac{\norm{\by}^2}{R_\by^2} + \frac{\norm{\tby}^2}{R_\tby^2} + \frac{\norm{\bz}}{R_\bz^2} + \frac{\norm{\tbz}^2}{R_\tbz^2}
\end{align*}
}
%for all $\bx, \tbx, \by, \tby, \bz, \tbz$ in their domains we have the following bound, which corresponds to the duality gap convergence of Extragradient algorithm:
According to the standard analysis of Mirror-Prox algorithm, the duality gap is bounded as follows:
\begin{equation}
    G_w(\bx_N, \tbx_N, \bz_N, \by, \tby, \tbz) - G_w(\bx, \tbx, \bz, \by_N, \tby_N, \tbz_N) \leq \ar{\frac{L_\zeta}{2N} \norm{\zeta - \zeta_0}^2},
\end{equation}
Substituting $\by = 0, \tby = 0, \tbz = 0, \bx = \bx_*, \tbx =\tbx_*, \bz = 0$ we get complexity estimate by function residual: 
\begin{equation}
    \sum_{k=1}^\ell f(x_N^k, \tx_N^k) - \sum_{k=1}^\ell f(x_*^k, \tx_*) \leq \frac{3L_\zeta}{N}.
\end{equation}
\ar{Analogously, we obtain bounds for affine constraints and consensus constraints
\begin{align*}
    \norm{Ax_N - b} + \norm{Cx_N - d} &\leq \frac{17\sqrt{2} L_\zeta}{N} \min_{k=1,\ldots,l}\sigma_{\min}^+({A^k}^\top, {C^k}^\top), \\
    \norm{\tilde\mA\tbx - \bb} + \norm{\tilde\mC\tbx - \bd} &\leq \frac{17\sqrt{2} L_\zeta}{N} \min_{k=1,\ldots,\ell} \sigma_{\min}^+({\tilde\mA}^\top, {\tilde\mC}^\top), \\
    \norm{\mW\by_N} &\leq \frac{17\sqrt{2}L_\zeta}{2 N} \lambda_{\min}^+(\mW), \\
    \norm{\tmW\tbx_N} &\leq \frac{17\sqrt{2}L_\zeta}{2 N} \lambda_{\min}^+(\tmW).
\end{align*}
}
% \begin{equation}
% \|Ax - b \| \leq ?, \|\tmA \tbx - \bb \| \leq ?, \|\tmW \tbx \| \leq ?
% \end{equation}

\begin{remark}\label{noneuclidian} 
In the problem formulation \eqref{eq:cpc:p} we can additionally assume that $x^k \in Q^k \subseteq \R^{n_k}$, $\tilde{x}\in \tilde{Q} \subseteq \R^{\tilde{n}}$, where $Q^k$ and $\tilde{Q}$ -- simple convex sets, i.e. simplex, ball, half plane e.t.c. In this case instead of decentralized Extragradient method for saddle point problem (Mirror Prox with euclidian prox-function) one should use general decentralized Mirror Prox algorithm \cite{rogozin2021decentralized}.
\end{remark}

\section{Numerical experiment}
For the purpose of numerical experiment data is taken from \cite{7428969}. Here 6 bus system contains 2 generators. DC optimal power flow  problem of the following form is considered:

\begin{subequations}
\begin{equation}
    \min_{\substack{p^G \in \sP \\  \theta \in \Theta ,}}
    \sum_{i\in G}c_i(p^G_i)
\end{equation}
\begin{equation}
    p^G_i-p^D_i = B_{ij}(\theta_i-\theta_j),
\end{equation}
\begin{equation}
   \left|(\theta_i - \theta_j) / X_{ij}\right| \leq F^{max}_{ij}.
\end{equation}
\end{subequations}

\textbf{Optimization variables:}
\begin{itemize}
    \item $p^G_i,\;i\in \{1,\ldots,l\} $ --- generator power output;
    % \item $p_{ij},\;i,j\in N$ --- power flow from bus $i$ to bus $j$;
    \item $\theta_i,\;i\in \{1,\ldots,l\}$ --- phase angle of the bus $i$.
\end{itemize}

\textbf{Parameters:}
\begin{itemize}
    \item $\sP = \prod_{i=1}^l \left[p^{G,\min}_i,~ p^{G,\max}_i \right]$ ~--- minimal and maximal generation. For nodes without generation $p^{G,\min}_i = p^{G,\max}_i = 0$.
    \item $\Theta =\prod_{i=1}^l \left[-\theta^{max}_i,~ \theta^{max}_i\right]$~--- maximal phase angle
    \item $p^D_i,\;i\in \{1,\ldots,l\}$ --- demand;
    \item $B_{ij} = B_{ji},\; i,j\in \{1,\ldots,l \}$ --- line susceptances. If no power line between nodes $i$ and $j$ then $B_{ij}=0$ else $B_{ij}>0$. $X_{ij} = -B_{ij},\; i,j\in \{1,\ldots,l\}$ are line reactances;
    % \item $p^{G,\max}_i,\;i\in G$ --- maximal generation;
    \item $F^{max}_{ij},\; i,j\in \{1,\ldots,l \}$ --- maximal power flow on the line $(i,j)$;
    % \item $\theta^{\max}_i,\;i\in N$ --- maximal phase angle.
\end{itemize}

\textbf{Cost functions:} 

$c_i(\cdot),\;i\in \{1,\ldots,l \}$ --- convex sufficiently smooth functions, representing the cost of operating a generator at given power.
\medskip

\begin{figure}
     \centering
     \begin{subfigure}[b]{0.45\textwidth}
         \centering
         \includegraphics[width=\textwidth]{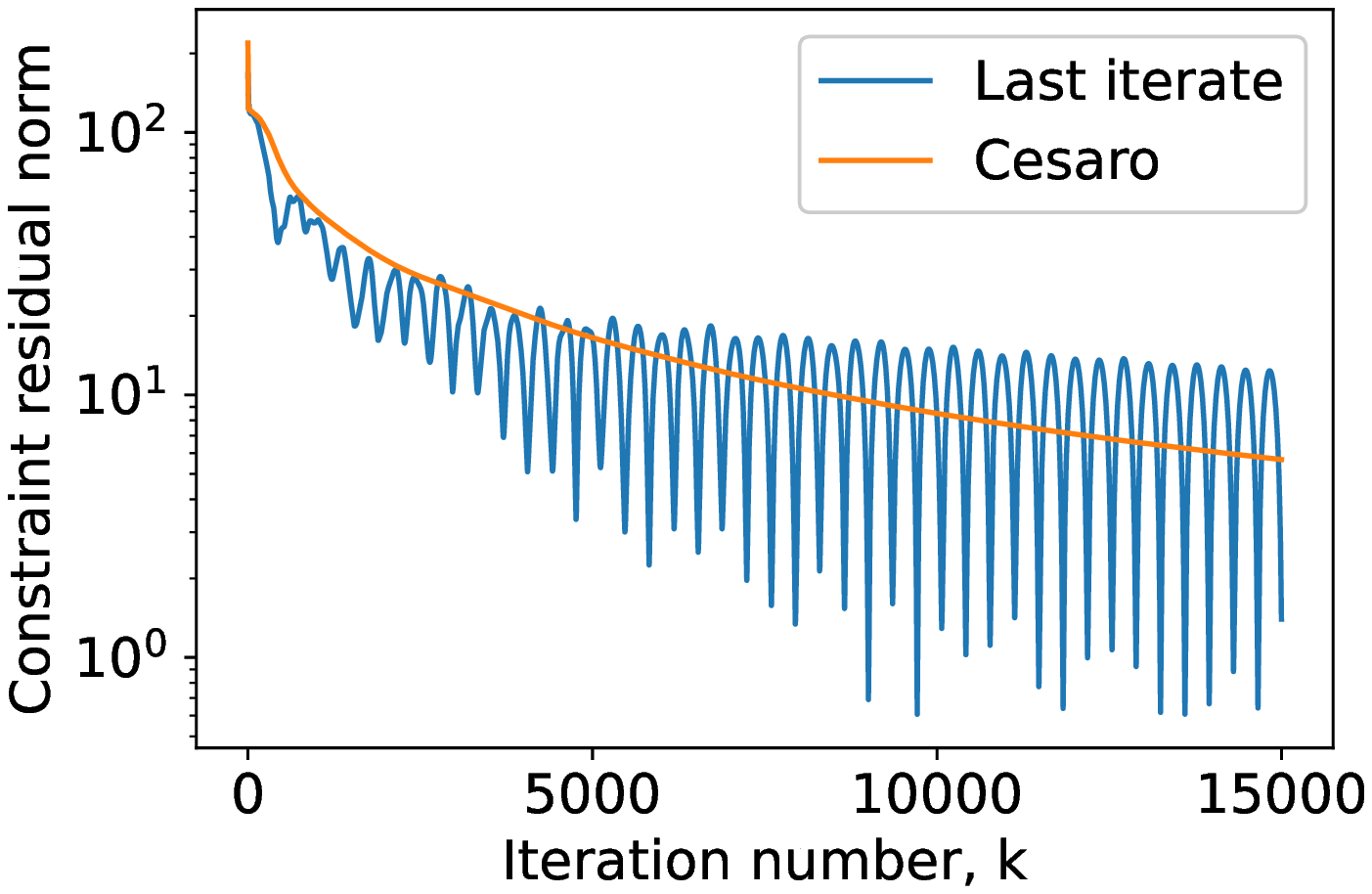}
        %  \caption{}
        %  \label{fig:no-cons:Ax-b}
     \end{subfigure}
     \begin{subfigure}[b]{0.45\textwidth}
         \centering
         \includegraphics[width=\textwidth]{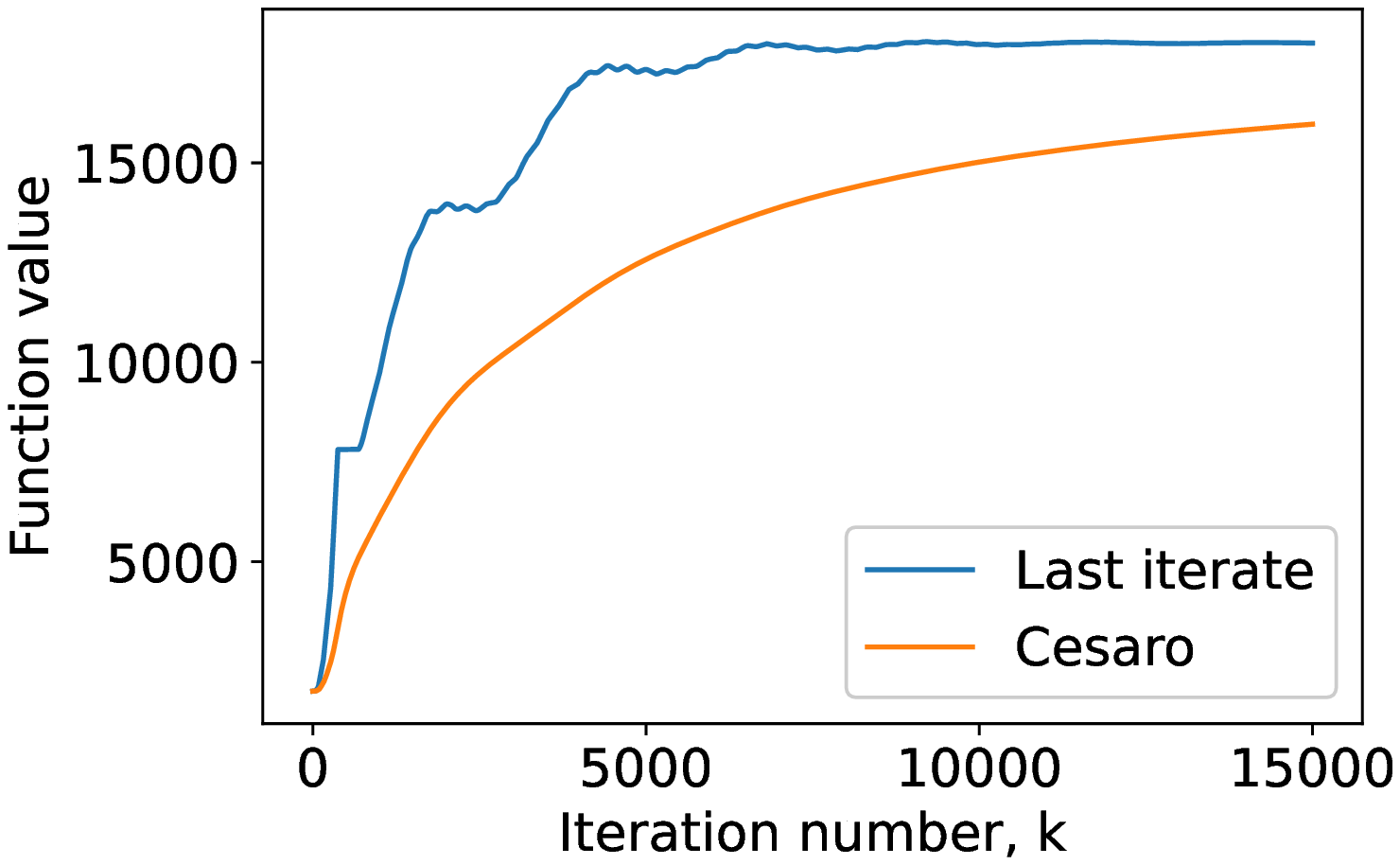}
        % \caption{}
        %  \label{fig:no-cons:f}
     \end{subfigure}
     \caption{Results of the numerical experiment for DC optimal power flow problem on 6-bus system \cite{7428969}}
\label{fig:no-cons}
\end{figure}

% {\color{red}Write equations that are in the program}.
The obtained results are consistent with the results in \cite{7428969}: generation is equal to 110 MW and 200 MW for the 1-st and 2-nd generators respectively. The results of numerical experiment are given in Fig. \ref{fig:no-cons}. 
Here the plots of function value and constraint residual convergence.

%%%%%%%%%%%%%%%%%%

% \include{samplepaper}

\bibliographystyle{splncs04}
\bibliography{bib}
\end{document}